\documentclass{llncs}
\usepackage{amsmath, amssymb, hyperref, color}
\usepackage{graphicx}

\newcommand{\RR}{\mathbb{R}}
\newcommand{\QQ}{\mathbb{Q}}

\begin{document}

\title{Maximum Likelihood Estimates for \\ Gaussian Mixtures Are Transcendental}

\author{
Carlos Am\'endola \inst{1} \and 
Mathias Drton \inst{2} \and Bernd Sturmfels \inst{1,3}}

\institute{Technische Universit\"{a}t, Berlin, Germany \and University of Washington, Seattle, 
USA \and University of California, Berkeley, USA}

\maketitle
 
\begin{abstract} 
Gaussian mixture models are central to classical statistics,
widely used in the information sciences, and
have a rich mathematical structure.  We examine
their maximum likelihood estimates  through
the lens of algebraic statistics. 
The MLE is not an algebraic function of the data, so 
there is no notion of ML degree for these models.
The critical points of the likelihood function are
transcendental, and there is no bound on their number, even 
for mixtures of 
two univariate Gaussians.

\end{abstract}

\section{Introduction}

The primary  purpose of this paper is to demonstrate
the result stated in the~title:

\begin{theorem}
\label{thm:main}
The maximum likelihood estimators of Gaussian mixture models are
transcendental functions.  More precisely, there exist rational samples
$x_1,x_2,\ldots,x_N$ in $\QQ^n$ whose maximum likelihood parameters for the
mixture of two $n$-dimensional Gaussians are not algebraic numbers over $\QQ$.
\end{theorem}

The principle of maximum likelihood (ML) is central to statistical
inference.  Most implementations of ML estimation employ iterative
hill-climbing methods, such as expectation maximization (EM).  These
can rarely certify that a globally optimal solution has been reached.
An alternative paradigm, advanced by algebraic statistics~\cite{DSS},
is to find the ML estimator (MLE) by solving the likelihood
equations. This is only feasible for small models, but it has the
benefit of being exact and certifiable.  
An important notion in this approach is
the {\em ML degree}, which is defined as the algebraic degree of the
MLE as a function of the data. This rests on the premise that the
likelihood equations are given by polynomials.

Many models used in practice, such as exponential families for
discrete or Gaussian observations, can be represented by polynomials.
Hence, they have an ML degree that serves as an upper bound for the
number of isolated local maxima of the likelihood function,
independently of the sample size and the data.  The ML degree is an
intrinsic invariant of a statistical model, with interesting geometric
and topological properties \cite{HS}.  The notion has proven useful
for characterization of when MLEs admit a `closed form'
\cite{SC}. When the ML degree is moderate, these exact tools are
guaranteed to find the optimal solution to the ML problem
\cite{BHR,GDP}.

However, the ML degree of a statistical model is only defined when the MLE is an algebraic
function of the data. Theorem \ref{thm:main} means that there is no ML
degree for Gaussian mixtures. It also highlights a fundamental
difference between likelihood inference and the method of moments
\cite{BS,GHK,Pea}. The latter is a computational paradigm within
algebraic geometry, that is, it is based on solving polynomial
equations.  ML estimation being transcendental means that likelihood
inference in Gaussian mixtures is outside the scope of algebraic
geometry.

The proof of Theorem \ref{thm:main} will appear in Section~\ref{sec2}.
In Section \ref{sec3} we shed further light on the transcendental
nature of Gaussian mixture models. 
We focus on mixtures of two univariate Gaussians, the model given in
(\ref{eq:pdf-uni}) below, and we present a family of data points on
the real line such that the number of critical points of the
corresponding log-likelihood function (\ref{eq:loglik1}) exceeds any
bound.

While the MLE for Gaussian mixtures is transcendental, this does not
mean that exact methods are not available.  Quite to the
contrary. Work of Yap and his collaborators in computational geometry
\cite{CCKPY,CPPY} convincingly demonstrates this. Using root bounds
from transcendental number theory, they provide certified answers to
geometric optimization problems whose solutions are known to be
transcendental.  Theorem \ref{thm:main} opens up the possibility of
transferring these techniques to statistical inference.  In our view,
Gaussian mixtures are an excellent domain of application for certified
computation in numerical analytic geometry.

\section{Reaching Transcendence}
\label{sec2}

Transcendental number theory \cite{Bak,Gel} is a field that
furnishes tools for deciding whether a given real number $\tau$
is a root of a nonzero polynomial in $\QQ[t]$.
If this holds then $\tau$ is algebraic; otherwise $\tau$ is transcendental.
For instance,  $\sqrt{2} + \sqrt{7} = 4.059964873...$ is algebraic, 
and so are the parameter estimates computed by Pearson in 
his 1894 study of crab data \cite{Pea}.
By contrast,  the famous constants $\pi = 3.141592653...$
and $e = 2.718281828...$ are  transcendental.
Our proof will be based on the following classical result.
A textbook reference is \cite[Theorem 1.4]{Bak}:

\begin{theorem}[Lindemann-Weierstrass]
If $u_1,\ldots,u_r$ are distinct algebraic numbers 
then $e^{u_1},\ldots,e^{u_r}$ are
linearly independent over the algebraic numbers.
\end{theorem}

For now, consider the case of $n=1$,
that is, mixtures of two univariate Gaussians.  We
allow mixtures with arbitrary means and variances.  Our model then
consists of all probability distributions on the real line $\RR$
with density
\begin{equation}
  \label{eq:pdf-uni}
  f_{\alpha,\mu,\sigma}(x)  = 
\frac{1}{\sqrt{2 \pi}} \cdot 
 \biggl[ \frac{\alpha}{\sigma_1} \,\exp \Bigl(-\frac{(x-\mu_1)^2 }{2\sigma_1^2} \Bigr)  \,+\,
\frac{1-\alpha}{\sigma_2} \,\exp \Bigl(-\frac{(x-\mu_2)^2 }{2\sigma_2^2} \Bigr)  
\biggr] \enspace .
\end{equation}
%
It has five unknown parameters, namely, the means
$\mu_1,\mu_2\in\mathbb{R}$, the standard deviations
$\sigma_1,\sigma_2>0$, and the mixture weight $\alpha\in[0,1]$.  The
aim 
is to estimate the five model parameters from
a collection of data points
$x_1,x_2,\ldots,x_N \in \RR$.  

The {\em log-likelihood function} of the model (\ref{eq:pdf-uni}) is
\begin{equation}
  \label{eq:loglik1}
 \ell(\alpha,\mu_1,\mu_2,\sigma_1,\sigma_2) \,\, = \,\,
\sum_{i=1}^N \log f_{\alpha,\mu,\sigma}(x_i) \enspace .
\end{equation}
This is a function of the five parameters, while
$x_1,\ldots,x_N$ are fixed constants.

The principle of maximum likelihood suggests to find estimates by
maximizing the function $\ell$ over the five-dimensional parameter
space $\Theta = [0,1] \times \RR^2 \times \RR^2_{>0}$.   

\begin{remark} 
\label{rmk:notbounded} \rm
The log-likelihood function $\ell$ in~(\ref{eq:loglik1}) is never bounded above.
To see this, we argue as in~\cite[Section~9.2.1]{Bishop}.
Set $N=2$,
fix arbitrary values
$\alpha_0\in[0,1]$, $\mu_{20}\in\mathbb{R}$ and $\sigma_{20}>0$, and
match the first mean to the first data point $\mu_1=x_1$.  The remaining
function of one unknown $\sigma_1$ equals
\[
  \ell(\alpha_0,x_1,\mu_{20},\sigma_1,\sigma_{20}) 
  \, \ge \, \log\biggl[\frac{\alpha_0}{\sigma_1} \,+\,
\frac{1-\alpha_0}{\sigma_{20}} \,\exp \Bigl(-\frac{(x_1-\mu_{20})^2
}{2\sigma_{20}^2}\Bigr)
\biggr]  \,+\, \text{const} \enspace .
\]
The lower bound tends to $\infty$ as $\sigma_1\to 0$.
\end{remark}

Remark \ref{rmk:notbounded} means that 
there is no global solution to the MLE problem. 
This is remedied by restricting to a subset of the parameter space
$\Theta$.  In practice, maximum likelihood for Gaussian mixtures means
computing local maxima of the function~$\ell$. These are found
numerically by a hill climbing method, such as the EM algorithm, with
particular choices of starting values. See Section~\ref{sec3}.
This method is implemented, for instance, in the {\tt R} package MCLUST
\cite{mclust}. In order for Theorem \ref{thm:main} to cover such local
maxima, we prove the following statement:  
\begin{quote}
  {\em There exist samples $x_1,\ldots,x_N\in\QQ$ such that every
    non-trivial critical point
    $(\hat \alpha, \hat \mu_1,\hat \mu_2,\hat \sigma_1,\hat \sigma_2)$
    of the log-likelihood function $\ell$ in the domain $\Theta$ has
    at least one transcendental coordinate.}
\end{quote}
Here, a critical point is {\em non-trivial } if it yields an honest
mixture, i.e.~a distribution that is not Gaussian.
By the identifiability results of \cite{teicher:1963}, this happens
if and only if the estimate
$(\hat \alpha, \hat \mu_1,\hat \mu_2,\hat \sigma_1,\hat \sigma_2)$ 
satisfies $0<\hat\alpha<1$ and
$(\hat\mu_1,\hat\sigma_1)\not=(\hat\mu_2,\hat\sigma_2)$.

\begin{remark} 
  The log-likelihood function always has some algebraic critical points, for any
  $x_1,\dots,x_N\in\QQ$.  Indeed, if we define the empirical
  mean and variance as
  \begin{align*}
    \bar x &=\frac{1}{N} \sum_{i=1}^N x_i \, ,
    & s^2 &= \frac{1}{N}\sum_{i=1}^N (x_i-\bar x)^2 \enspace ,
  \end{align*}
  then any point
  $(\hat \alpha, \hat \mu_1,\hat \mu_2,\hat \sigma_1,\hat \sigma_2)$
  with $\hat\mu_1=\hat\mu_2=\bar x$ and $\hat\sigma_1=\hat\sigma_2=s$
  is critical.  This gives a Gaussian distribution
  with mean $\bar x$ and variance $s^2$, so it is trivial.
\end{remark}

\begin{proof}[of Theorem \ref{thm:main}]
  First, we treat the univariate case. 
  Consider the partial derivative  of (\ref{eq:loglik1})
  with respect to the mixture weight $\alpha$:
  \begin{small}
  \begin{align}
    \label{eq:partial-alpha}
    \frac{\partial\ell}{\partial\alpha}
    &\;=\; \sum_{i=1}^N \frac{1}{f_{\alpha,\mu,\sigma}(x_i)}\cdot
      \frac{1}{\sqrt{2\pi}}\left[
      \frac{1}{\sigma_1}\exp \Bigl(-\frac{(x_i-\mu_1)^2
      }{2\sigma_1^2}\Bigr)\;-\; 
\frac{1}{\sigma_2} \exp \Bigl(-\frac{(x_i-\mu_2)^2 }{2\sigma_2^2}
      \Bigr) \right].                                                 
  \end{align}
  \end{small}
  Clearing the common denominator
  \[
  \sqrt{2\pi}\cdot\prod_{i=1}^N f_{\alpha,\mu,\sigma}(x_i) \enspace ,
  \]
  we see that $\partial\ell/\partial\alpha=0$ if and only if
  \begin{align}
    \nonumber
    \sum_{i=1}^N &\left[
      \frac{1}{\sigma_1}\,\exp \Bigl(-\frac{(x_i-\mu_1)^2
      }{2\sigma_1^2}\Bigr)\;-\; 
\frac{1}{\sigma_2} \,\exp \Bigl(-\frac{(x_i-\mu_2)^2 }{2\sigma_2^2}
      \Bigr) \right]\\
    \label{eq:partial-alpha-cleared}
    & \times\prod_{j\not=i} \biggl[ \frac{\alpha}{\sigma_1} \,\exp
      \Bigl(-\frac{(x_j-\mu_1)^2 }{2\sigma_1^2} \Bigr)  \,+\, 
\frac{1-\alpha}{\sigma_2} \,\exp \Bigl(-\frac{(x_j-\mu_2)^2
      }{2\sigma_2^2} \Bigr)   
\biggr]\;=\; 0 \enspace .
  \end{align}
  Letting $\alpha_1=\alpha$ and $\alpha_2=1-\alpha$, we may rewrite
  the left-hand side of~(\ref{eq:partial-alpha-cleared}) as
  \begin{small}
  \begin{align}
    \label{eq:partial-alpha-cleared2}
    \sum_{i=1}^N &\left[  \sum_{k_i=1}^2 \frac{(-1)^{k_i-1}}{\sigma_{k_i}}\,\exp \Bigl(-\frac{(x_i-\mu_{k_i})^2
                   }{2\sigma_{k_i}^2}\Bigr)\right]
                   \prod_{j\not=i} \left[ \sum_{k_j=1}^2
                   \frac{\alpha_{k_j}}{\sigma_{k_j}} \,\exp 
      \Bigl(-\frac{(x_j-\mu_{k_j})^2 }{2\sigma_{k_j}^2} \Bigr) 
\right]. 
  \end{align}
  \end{small}
  We expand the products, collect terms, and set
    $N_i(k)=\vert \{j:k_j=i\} \vert$. With this,   the partial derivative
  $\partial\ell/\partial\alpha$ is zero if and only if the following vanishes:
  \begin{small}
  \begin{align*}
    \nonumber
    &\sum_{i=1}^N \sum_{k\in\{1,2\}^N}\!\!\!
    \exp \! \left( \! -\sum_{j=1}^N \frac{(x_j-\mu_{k_j})^2
      }{2\sigma_{k_j}^2} \right)\! (-1)^{k_i-1}\alpha^{|\{j\not=i:
    k_j=1\}|}(1{-}\alpha)^{|\{j\not=i: k_j=2\}|} \left(\prod_{j=1}^N
    \frac{1}{\sigma_{k_j}}  \right)\\
    \nonumber
    &\! =\!\!\!  \sum_{k\in\{1,2\}^N}\!\!\!\!
    \exp \! \left(\! -\sum_{j=1}^N \!\frac{(x_j-\mu_{k_j})^2
      }{2\sigma_{k_j}^2} \right)  \!\!\left(\prod_{j=1}^N
    \frac{1}{\sigma_{k_j}}  \right) \! \sum_{i=1}^N (-1)^{k_i-1}\alpha^{|\{j\not=i:
    k_j=1\}|}(1-\alpha)^{|\{j\not=i: k_j=2\}|}\\
    \nonumber
    &\!\! = \!\! \sum_{k\in\{1,2\}^N} \!\!
\!\!    \exp \!\left( -\sum_{j=1}^N \frac{(x_j-\mu_{k_j})^2
      }{2\sigma_{k_j}^2} \right)\! \! \left(\prod_{j=1}^N
    \frac{1}{\sigma_{k_j}}  \right) \!
      \alpha^{N_1(k)-1}(1{-}\alpha)^{N_2(k)-1}
      \left[ { N_1(k)(1{-}\alpha) \atop +N_2(k)(-\alpha) } \right]\\
    &\! \! = \!\! \sum_{k\in\{1,2\}^N} \!\!\!\!
    \exp \! \left( \! -\sum_{j=1}^N \frac{(x_j-\mu_{k_j})^2
      }{2\sigma_{k_j}^2} \right) \!\! \left(\prod_{j=1}^N
    \frac{1}{\sigma_{k_j}}  \right)
      \alpha^{N_1(k)-1}(1-\alpha)^{N-N_1(k)-1}(N_1(k)-N\alpha) \enspace . \,\,\,
  \end{align*}
\end{small}
    
  Let  $(\hat \alpha, \hat \mu_1, \hat \mu_2, \hat \sigma_1, \hat   \sigma_2)$
  be a non-trivial isolated critical point of the likelihood function. This means that $0<\hat\alpha<1$ and
  $(\hat\mu_1,\hat\sigma_1)\not=(\hat\mu_2,\hat\sigma_2)$. 
This point depends continuously on the choice of the data 
    $x_1,x_2,\ldots,x_N$. By moving the vector with these coordinates
    along a general line in $\RR^N$, the mixture parameter $\hat \alpha$
    moves continuously in the critical equation $\partial\ell/\partial\alpha = 0$ above. By the Implicit Function Theorem, it takes on all values 
    in some open interval  of $\RR$, and we can thus choose our data points $x_i$
    general enough so that
     $\hat\alpha$ is not an integer multiple
  of $1/N$. We can further ensure that 
  the last sum above is a 
  $\mathbb{Q}(\alpha)$-linear
  combination of exponentials with nonzero coefficients. 
  
  Suppose that $(\hat \alpha, \hat \mu_1, \hat \mu_2, \hat \sigma_1, \hat
  \sigma_2)$ is algebraic. The
  Lindemann-Weierstrass Theorem implies that the 
  arguments of ${\rm exp}$ 
     are all the same. 
    Then the   $2^N$ numbers 
      \[
  \sum_{j=1}^N \frac{(x_j-\hat \mu_{k_j})^2
      }{2{\hat \sigma_{k_j}}^2} \, ,\qquad k\in\{1,2\}^N,
  \]
  are all identical.  However, for $N \geq 3$, and for general choice
  of data $x_1,\dots,x_N$ as above, this can only happen if
  $(\hat\mu_1,\hat\sigma_1)=(\hat\mu_2,\hat\sigma_2)$.  This
  contradicts our hypothesis that the critical point is non-trivial.
  We conclude that all non-trivial critical points of the
  log-likelihood function (\ref{eq:loglik1}) are transcendental. 
  \smallskip

  In the multivariate case, the model parameters comprise the mixture
  weight $\alpha\in[0,1]$,
  mean vectors $\mu_1,\mu_2\in\RR^n$
  and positive definite covariance matrices
  $\Sigma_1,\Sigma_2\in\RR^n$.
  Arguing as above, if a non-trivial critical
  $(\hat \alpha, \hat \mu_1, \hat \mu_2, \hat \Sigma_1, \hat
  \Sigma_2)$ is algebraic, then the Lindemann-Weierstrass Theorem
  implies that the numbers
  \[
    \sum_{j=1}^N (x_j-\hat \mu_{k_j})^T \hat
    \Sigma_{k_j}^{-1}(x_j-\hat \mu_{k_j}) \, ,\qquad k\in\{1,2\}^N \, ,
  \]
  are all identical.  For $N$
  sufficiently large and a general choice of $x_1,\dots,x_N$
  in $\RR^n$,
  the  $2^N$ numbers are identical only if
  $(\hat\mu_1,\hat\Sigma_1)=(\hat\mu_2,\hat\Sigma_2)$.
  Again, this constitutes a contradiction to the hypothesis that
  $(\hat \alpha, \hat \mu_1, \hat \mu_2, \hat \Sigma_1, \hat
  \Sigma_2)$ is non-trivial.
\qed
\end{proof}

Many variations and specializations of the Gaussian mixture model
are used in applications.  In the case $n=1$, the
variances are sometimes assumed equal, so
$\sigma_1=\sigma_2$ for the above two-mixture.  This avoids the issue
of an unbounded likelihood function (as long as $N\ge 3$).  Our proof
of Theorem~\ref{thm:main}  applies to this setting.  In higher
dimensions $(n\ge 2$), the covariance matrices are sometimes assumed
arbitrary and distinct, sometimes arbitrary and equal, but often also
have special structure such as being diagonal.
Various default choices are discussed in the paper \cite{mclust} that 
introduces the  {\tt R} package MCLUST.  
Our results imply that
maximum likelihood estimation is transcendental
for all these MCLUST models. 

\begin{figure}[t]
  \centering
 \includegraphics[scale=0.45]{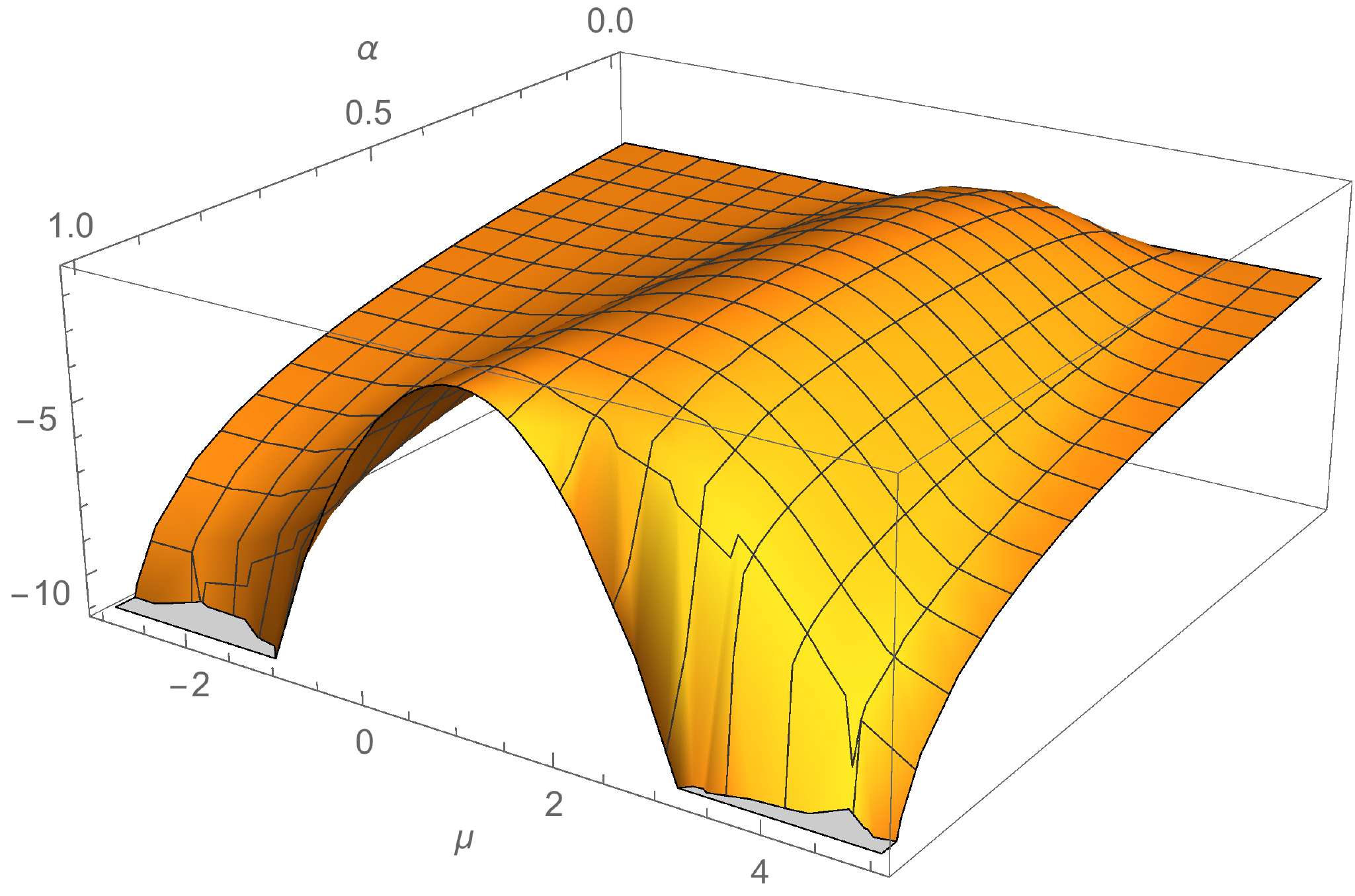}
\vspace{-0.1in}
  \caption{Graph of the log-likelihood function 
         for two data points $x_1=0$ and $x_2=2$.}
  \label{fig:simple}
\end{figure}

\begin{example} \rm We illustrate
  Theorem~\ref{thm:main} for a 
specialization of
  (\ref{eq:pdf-uni}) obtained by fixing three 
 parameters:
     $\mu_2=0$ and $\sigma_1=\sigma_2=1/\sqrt{2}$.  The
  remaining two free parameters are $\alpha$ and $\mu=\mu_1$.  We take only
  $N=2$ data points, namely $x_1 = 0$ and $x_2 = x>0$.  
  Omitting
  an additive constant, our log-likelihood function equals
\begin{equation}
  \label{eq:simplest-loglik-0x}
  \ell(\alpha,\mu) \,\,\, = \,\, \,
  {\rm log}\bigl(  \alpha \cdot e^{ - \mu^2 }+ (1-\alpha)   \bigr) \,\, + \,\,
  {\rm log}\bigl(  \alpha \cdot e^{ - (\mu-x)^2 }+ (1-\alpha) \cdot e^{-x^2}  \bigr)\enspace .
\end{equation}

For a concrete example take $x = 2$. The graph of 
(\ref{eq:simplest-loglik-0x}) for this choice is shown in
Figure~\ref{fig:simple}.
By maximizing $\ell(\alpha,\mu)$ numerically, we find
the parameter estimates
\begin{equation}
\label{eq:two numbers}
\hat \alpha =  0.50173262959803874...
\quad \hbox{and} \quad 
\hat \mu \, = \, 1.95742494230308167...
\end{equation}
Our technique can be applied to prove that $\hat \alpha$ 
and $\hat \mu$ are transcendental over $\QQ$. We illustrate it for $\hat \mu$. 

For any $x\in\RR$, the function $\ell(\alpha,\mu)$ is bounded from
above and achieves its maximum on $[0,1]\times \mathbb{R}$.  
If $x>0$ is large, then any global maximum
$(\hat\alpha,\hat\mu)$ of $\ell$ is in the interior of
$[0,1]\times \mathbb{R}$ and satisfies $0<\hat\mu\le x$.  According to
a {\tt Mathematica} computation, the choice $x\ge 1.56125...$ suffices for
this.  Assume that this holds.  Setting the two partial derivatives
equal to zero and eliminating the unknown $\alpha$ in a further {\tt
  Mathematica} computation, the critical equation for $\mu$ is found
to be
\begin{equation}
  \label{eq:transeqn}
   (x-\mu)e^{\mu^2}  \,-\,x \,+\, \mu e^{-\mu (2x-\mu)} \,\;=\,\;0 \enspace .
\end{equation}
Suppose for contradiction that both $x$ and $\hat\mu$ are algebraic
numbers over $\QQ$.  Since $0<\hat\mu\le x$, we have 
$-\hat\mu (2x-\hat\mu) < 0 < \hat\mu^2$.  Hence $\,u_1 = \hat\mu^2$,
$\,u_2 = 0\,$ and $\,u_3 = -\hat\mu (2x-\hat\mu)\,$ are distinct
algebraic numbers.  The Lindemann-Weierstrass Theorem implies that
$e^{u_1},e^{u_2}$ and $e^{u_3}$ are linearly independent over the
field of algebraic numbers.  However, from (\ref{eq:transeqn}) we know
that
\[
(x-\hat\mu) \cdot e^{u_1} \,-\,x \cdot e^{u_2} 
\,+\, \hat\mu \cdot e^{u_3} \,\,\;=\,\,\;0 \enspace .
\]
This is a contradiction. We
conclude that the number $\hat\mu$ is transcendental over~$\QQ$.
\end{example}

\section{Many Critical Points}
\label{sec3}

Theorem \ref{thm:main} shows that Gaussian mixtures
do not admit an  ML degree. This raises the question of how
to find any bound for the number of critical points.

\begin{problem}
Does there exist a universal bound on the number of non-trivial critical points for the log-likelihood function of the mixture
of two univariate Gaussians? Or, can we find
a sequence of samples on the real line such that
the number of non-trivial critical points increases beyond any bound?
\end{problem}

We shall resolve this problem by answering the
second question affirmatively.  The idea behind our solution is to choose a sample consisting of
  many well-separated clusters of
 size $2$.  Then each cluster gives rise to a distinct
 non-trivial critical point
 $(\hat{\alpha},\hat{\mu}_1,\hat{\mu}_2,\hat{\sigma}_1,\hat{\sigma}_2)$
 of the log-likelihood function $\ell$ from (\ref{eq:loglik1}).  
  We
propose one particular choice of data, but many others
would work too.

 \begin{theorem}
 \label{thm:manyhills}
 Fix sample size $N=2K$ for $K \geq 2$, and take the ordered sample
 $(x_1,\ldots,x_{2K})=(1,1.2,\,2,2.2,\,\ldots,K,K{+}0.2)$.  Then, for
 each $k \in \{1,\ldots,K\}$, the log-likelihood function $\ell$ from
{\rm (\ref{eq:loglik1})} has a non-trivial critical point with
 $ k < \hat \mu_1 < k +0.2$.  Hence, there are at least $K$
 non-trivial critical points.
\end{theorem}

Before turning to the proof, we offer a numerical illustration.

\begin{table}[t]
  \caption{\label{tab:K7} Seven critical points of the log-likelihood
    function in Theorem~\ref{thm:manyhills} with $K=7$.}
\begin{center} \begin{tabular}{ | l | l | l | l | l | l | p{3.1cm} |} \hline  $k$ & $\alpha$ & $\mu_1$ & $\mu_2$ & $\sigma_1$ & $\sigma_2$ & log-likelihood\\ 
\hline 1 & 0.1311958 & 1.098998 & 4.553174 & 0.09999497 & 1.746049 & -27.2918782147578 \\
\hline 2 & 0.1032031 & 2.097836 & 4.330408 & 0.09997658 & 1.988948 & -28.6397463805501 \\
\hline 3 &  0.07883084 & 3.097929 & 4.185754 & 0.09997856 & 2.06374 & -29.1550277534757 \\
\hline 4 & 0.06897294 & 4.1 & 4.1 &0.1 & 2.07517 & -29.2858981551065 \\
\hline 5 & 0.07883084 & 5.102071 & 4.014246 & 0.09997856 & 2.06374 & -29.1550277534757 \\
\hline 6 & 0.1032031 & 6.102164 & 3.869592 & 0.09997658 & 1.988948 & -28.6397463805501 \\
\hline 7 & 0.1311958 & 7.101002 & 3.646826 & 0.09999497 & 1.746049 & -27.2918782147578 \\
\hline \end{tabular} 
\end{center}
\end{table}

\begin{example} \rm For $K = 7$, we have $N=14$ data points in the
  interval $[1,7.2]$.
  Running the EM algorithm (as explained in the proof of Theorem~\ref{thm:manyhills}
  below)  yields
  the non-trivial critical points 
  reported 
  in Table~\ref{tab:K7}.  Their $\mu_1$ coordinates are seen
  to be close to the cluster midpoints $k+0.1$ for all $k$.
    The observed symmetry under reversing the order of the
  rows also holds for all larger $K$.
\end{example}

Our proof of Theorem \ref{thm:manyhills} will be based on the {\em EM
  algorithm}.   We first recall this  algorithm.  Let
$f_{\alpha,\mu,\sigma}$ be the mixture density
from~(\ref{eq:pdf-uni}), and let
\[
f_j(x)\,\,=\,\,\frac{1}{\sqrt{2 \pi}\, \sigma_j}  \,\exp
\Bigl(-\frac{(x-\mu_j)^2 }{2\sigma_j^2} \Bigr), \quad j=1,2 \enspace ,
\]
be the two Gaussian component densities.  Define
\begin{equation} \label{probs}
\gamma_i=\frac{\alpha \cdot f_1(x_i)}{f_{\alpha,\mu,\sigma}(x_i)} \enspace ,
\end{equation}
which can be interpreted as the conditional probability that data
point $x_i$ belongs to the first mixture component.  Further, define
$N_1=\sum_{i=1}^N \gamma_i$ and $N_2=N-N_1$, which are expected
cluster sizes.  Following \cite[Section~9.2.2]{Bishop}, the likelihood
equations for our model can be written in the following fixed-point form:
\begin{align}
 \alpha &= \frac{N_1}{N} \, , \label{firstlik}\\
 \mu_1 &= \frac{1}{N_1}\sum_{i=1}^N \gamma_ix_i \, , \label{mu1} &
  \mu_2 &= \frac{1}{N_2}\sum_{i=1}^N (1-\gamma_i)x_i \, ,  \\
\sigma_1 &= \frac{1}{N_1}\sum_{i=1}^N
           \gamma_i(x_i-\mu_1)^2 \, , &
  \sigma_2 &= \frac{1}{N_2}\sum_{i=1}^N
             (1-\gamma_i)(x_i-\mu_2)^2 \enspace .\label{lastlik}
\end{align}
In the present context, the EM algorithm amounts to solving these
equations iteratively.
More precisely, consider any starting point
$(\alpha,\mu_1,\mu_2,\sigma_1,\sigma_2)$. Then the E-step
(``expectation'') computes the estimated frequencies $\gamma_i$ via
(\ref{probs}).  In the subsequent M-step (``maximization''), one
obtains a new parameter vector
$(\alpha,\mu_1,\mu_2,\sigma_1,\sigma_2)$ by evaluating the right-hand
sides of the equations (\ref{firstlik})-(\ref{lastlik}).  The two
steps are repeated until a fixed point is reached, up to the desired
numerical accuracy.  The updates never decrease the log-likelihood.
For our problem it can be shown that the algorithm will converge to a
critical point; see e.g.~\cite{RW}.
  
\begin{proof}[of Theorem \ref{thm:manyhills}]
  Fix  $k\in\{1,\dots,K\}$.  We
  choose starting parameter values to suggest that the pair
  $(x_{2k-1},x_{2k})=(k,k+0.2)$ belongs to the first mixture
  component, while the rest of the sample belongs to the second.
  Explicitly, we set
  \begin{align*}
    \alpha &= \frac{2}{N} = \frac{1}{K} \, ,\\
    \mu_1 &= k+0.1 \, ,  &
                       \mu_2 &= \frac{K^2 + 1.2K -2k -0.2}{2(K-1)} \, , \\
    \sigma_1 &= 0.1 \, , &
                       \sigma_2 &= \frac{\sqrt{\frac{1}{12}K^4 -\frac{1}{3}K^3 + (k - \frac{43}{75})K^2 -(k^2 -k +\frac{14}{75})K+0.01}}{K-1} \enspace .
  \end{align*}
  We shall argue that, when running the EM algorithm, the parameters
  will always stay close to these starting values.
  Specifically, we claim that throughout all EM iterations, the
  parameter values satisfy the inequalities
  \begin{align}
    \label{eq:box1}
    &\frac{1}{4K} \leq \alpha \leq \frac{1}{K} \, , &\
    &0.09 \leq \mu_1-k \leq 0.11 \, , &
    &0.099 \leq \sigma_1 \leq 0.105 \, , 
  \end{align}
  \vspace{-0.21in}
  \begin{align}
    \label{eq:box2}
    \frac{K}{2}+0.1 \,\,\leq \,\,&\mu_2 \leq \frac{K}{2}+1.1 \, ,\\
    \label{eq:box3}
    \sqrt{\frac{K^2}{12}-\frac{K}{6}+0.01} \,\,\leq \,\,
      &\sigma_2 \leq \sqrt{\frac{K^2}{12}+\frac{K}{12}+0.01} \enspace .
  \end{align}
  The starting values proposed above obviously satisfy the
  inequalities in~(\ref{eq:box1}), and it is not difficult to check
  that~(\ref{eq:box2}) and~(\ref{eq:box3}) are satisfied as well.  To
  prove the theorem, it remains to show
  that~(\ref{eq:box1})-(\ref{eq:box3}) continue to hold after an
  EM update.

  In the remainder, we assume that $K>22$.  For smaller values of $K$
  the claim of the theorem can be checked by running the EM algorithm.
  In particular, for $K>3$, the second standard deviation satisfies
  the simpler bounds
  \begin{equation}
    \label{eq:box4}
  \frac{K}{\sqrt{12}} - \frac{\sqrt{3}}{5} \leq \sigma_2 \leq
  \frac{K}{\sqrt{12}}+ \frac{\sqrt{3}}{12} \enspace .
  \end{equation}

  A key property is that the quantity $\gamma_i$, computed in the
  E-step, is always very close to zero for $i \neq 2k-1,2k$.  To see
  why, rewrite~(\ref{probs}) as
  \[
  \gamma_i=
  \frac{1}{1+\frac{1-\alpha}{\alpha}\frac{f_2(x_i)}{f_1(x_i)}}=\frac{1}{1+\frac{1-\alpha}{\alpha}\frac{\sigma_1}{\sigma_2}\exp
    \left\lbrace \frac{1}{2}\left( (\frac{x_i - \mu_1}{\sigma_1})^2 -
        (\frac{x_i - \mu_2}{\sigma_2})^2 \right) \right\rbrace} \enspace .
  \] 
  Since $\alpha \leq 1/K$, we have $\frac{1-\alpha}{\alpha} \geq K-1$.
  On the other hand,
  $\frac{\sigma_1}{\sigma_2} \geq \frac{0.099}{K / \sqrt{12} +
    \sqrt{3} / 12}$.
  Using that $K>22$, their product is thus bounded below by 0.3209.
  Turning to the exponential term, the second inequality
  in~(\ref{eq:box1}) implies that $|x_i-\mu_1|\ge 0.89$ for $i = 2k-2$
  or $i=2k+1$, which index the data points closest to the $k$th pair.
  
  Using~(\ref{eq:box4}), we obtain 
  \[
  \left(\frac{x_i - \mu_1}{\sigma_1}\right)^2 - \left(\frac{x_i -
      \mu_2}{\sigma_2}\right)^2 
  \,\, \geq \,\, \left(\frac{0.89}{0.105}\right)^2 - \left(\frac{K/2+0.1}{K /
    \sqrt{12} - \sqrt{3} / 5}\right)^2 \,\, \geq \,\, 67.86 \enspace .
  \]
  From $e^{33.93} > 5.4 \cdot 10^{14}$, we deduce that
  $\gamma_i < 10^{-14}$.  The exponential term becomes only smaller as
  the considered data point $x_i$ move away from the $k$th pair. As
  $| i-(2k-1/2)|$ increases, $\gamma_i$ decreases and can be bounded
  above by a geometric progression starting at $10^{-14}$ and with
  ratio $10^{-54}$.  This makes $\gamma_i$ with $i \neq 2k, 2k-1$
  negligible. Indeed, from the limit of geometric series, we have
  \begin{align}
    s_1 & = \sum_{i\not=2k-1,2k} \gamma_i  \,\, <\,\,10^{-13} \enspace ,
  \end{align}
  and similarly, $s_2=\sum_{i\not=2k-1,2k} \gamma_i(x_i-k)$ satisfies
  \begin{equation}
  \label{series}
\!\!    | s_2| \,=\, \left| \gamma_{2k-2}(-0.8) + \gamma_{2k+1}(1) +
             \gamma_{2k-3}(-1) + \gamma_{2k+2}(1.2) + \ldots \right|
             \, <\, 10^{-13} \enspace .   
             \end{equation}
  The two sums $s_1$ and $s_2$ are relevant for the M-step.

  The probabilities $\gamma_{2k-1}$ and $\gamma_{2k}$  give
  the main contribution to the averages that are evaluated in the M-step.
  They satisfy $\,  0.2621 \leq \gamma_{2k-1}, \gamma_{2k} \leq 0.9219$.
  Moreover, we may show that the values of $\gamma_{2k-1}$ and
  $\gamma_{2k}$ are similar, namely:
  \begin{equation}
    \label{eq:bound:gamma:ratio}
  0.8298 \leq \frac{\gamma_{2k-1}}{\gamma_{2k}}\leq 1.2213 \enspace ,
  \end{equation}
  which we prove by writing
  \[
  \frac{\gamma_{2k-1}}{\gamma_{2k}} = \frac{1+y \exp(z/2)}{1+y} \enspace ,
  \]
  and using $K>22$ to bound
  \begin{eqnarray*}
    y &\,=\,& \frac{1-\alpha}{\alpha}\frac{\sigma_1}{\sigma_2}\exp
              \left\lbrace \frac{1}{2}\left( \left(\frac{k -
              \mu_1}{\sigma_1}\right)^2 - \left(\frac{k -
              \mu_2}{\sigma_2}\right)^2 \right) \right\rbrace,  \\ 
    z &\,=\,& \frac{0.4(k-\mu_1)+0.04}{\sigma_1^2} -
          \frac{0.4(k-\mu_2)+0.04}{\sigma_2^2} \enspace . 
  \end{eqnarray*}

  Bringing it all together, we have
  $$
  \mu_1 = \frac{1}{N_1}\sum_{i=1}^N \gamma_ix_i = k+\frac{0.2
    \gamma_{2k}+ s_2}{\gamma_{2k-1} + \gamma_{2k} + s_1} \enspace . 
  $$
  Using  $\gamma_{2k}+\gamma_{2k-1}>0.5$ and
  (\ref{series}), as well as  the lower bound in (\ref{eq:bound:gamma:ratio}),
  we find
  $$  \mu_1 - k \leq \frac{0.2 \gamma_{2k}}{\gamma_{2k-1} +
    \gamma_{2k}} + \frac{s_2}{\gamma_{2k-1} + \gamma_{2k}} \leq
  \frac{0.2\gamma_{2k}}{0.8298\gamma_{2k}+\gamma_{2k}}+10^{-12} \leq
  0.11 \enspace .
  $$
  Using the upper bound in (\ref{eq:bound:gamma:ratio}), we also have
  $0.09 \leq \mu_1 - k$.  Hence, the second inequality in~(\ref{eq:box1}) holds.
 
  The inequalities for the other parameters are verified similarly.
  For instance,   $$
  \frac{1}{4K} < \frac{0.2621 + 0.2621}{2K} \leq  \frac{\gamma_{2k-1}
    + \gamma_{2k} + s_1}{2K} \leq \frac{0.9219 + 0.9219 +
    10^{-13}}{2K} < \frac{1}{K}
  $$
  holds for $\alpha =\frac{N_1}{N}$.
  Therefore, the first inequality in~(\ref{eq:box1})
  continues to be true.

  We conclude that running the EM algorithm from the chosen starting
  values yields a sequence of parameter vectors that satisfy
  the inequalities~(\ref{eq:box1})-(\ref{eq:box4}).  
  The
  sequence has at least one limit point, which must be a non-trivial
  critical point of the log-likelihood function.  Therefore, for every
  $k=1,\dots,K$, the log-likelihood function has a non-trivial
  critical point with $\mu_1\in(k,k+0.2)$. \qed
\end{proof}
	
\section{Conclusion}	

We showed that the maximum likelihood estimator (MLE) in Gaussian
mixture models is not an algebraic function of the data, and that
the log-likelihood function may have
arbitrarily many critical points.  Hence, in contrast to the models
studied so far in algebraic statistics \cite{BHR,DSS,GDP,SC}, there is
no notion of an ML degree for Gaussian mixtures.  However, certified
likelihood inference may still be possible, via transcendental root
separation bounds, as in \cite{CCKPY,CPPY}.
  
\begin{remark}
  \rm The {\em Cauchy-location model}, treated in \cite{Ree}, is an example
  where the ML estimation is algebraic but the ML degree, and also the
  maximum number of local maxima, depends on the sample size and
  increases beyond any bound.
\end{remark}

\begin{remark}
  \rm The ML estimation problem admits a population/infinite-sample
  version.   Here the maximization of the
  likelihood function is replaced
   by {\em minimization of the Kullback-Leibler
  divergence} between a given data-generating distribution and the
  distributions in the model.  The question of whether this
  population problem is subject to local but not global maxima was
  raised in \cite{Sre}---in the context of Gaussian mixtures with
  known and equal variances.  It is known that the Kullback-Leibler
  divergence for such Gaussian mixtures is not an analytic function
  \cite[\S7.8]{WatanabeBook}.  Readers of Japanese should be able to
  find details in \cite{WatanabeNonanalytic}.
\end{remark}

As previously mentioned, Theorem \ref{thm:main} shows that likelihood
inference is in a fundamental way more complicated than the 
classical {\em method of moments} \cite{Pea}.
The latter involves only the solution of
polynomial equation systems. This was recognized also 
in the computer science literature on learning Gaussian mixtures
\cite{BS,GHK,MV}, where most of the recent progress
is based on variants of the method of moments
rather than likelihood inference. 
We refer to \cite{AFS} for a study of the method
of moments from an algebraic perspective.
Section 3 in that paper illustrates
the behavior of Pearson's method for the sample 
used in Theorem~\ref{thm:manyhills}.

\bigskip
\medskip

\noindent
{\bf Acknowledgements.}\smallskip \\
CA and BS were supported by the Einstein Foundation Berlin.
MD and BS also thank the US National Science Foundation 
(DMS-1305154 and DMS-1419018).

\end{document}